\documentclass{amsart}

\usepackage[mathscr]{eucal}
\usepackage{amsfonts}
\usepackage{amsmath}
\usepackage{amsthm}
\usepackage{amssymb}
\usepackage{latexsym}
\usepackage{epsfig}
\usepackage{hhline}
\input xy
\xyoption{all}

\newtheorem{theorem}{Theorem}[section]
\newtheorem{lemma}[theorem]{Lemma}
\newtheorem{proposition}[theorem]{Proposition}

\theoremstyle{definition}
\newtheorem{definition}[theorem]{Definition}
\newtheorem{claim}[theorem]{Claim}

\newtheorem{example}[theorem]{Example}

\newtheorem{remark}[theorem]{Remark}

\newcommand{\Ker}{\text{Ker\,}}
\newcommand{\rank}{\text{rank\,}}

\newcommand{\Hom}{\operatorname{Hom}}

\newcommand{\Id}{\text{Id}}

\newcommand{\Rep}{\operatorname{Rep}}

\newcommand{\tr}{\text{tr}\,}

\newcommand{\ben}{\begin{enumerate}}
\newcommand{\een}{\end{enumerate}}
\newcommand{\beq}{\begin{equation}}
\newcommand{\eeq}{\end{equation}}

\newcommand{\CC}{{\mathbb{C}}}
\newcommand{\ZZ}{{\mathbb{Z}}}

\newcommand{\sgn}{{\sf {sgn}}}

\begin{document}
\title{Twisted symplectic reflection algebras}

\begin{abstract}
In this paper we introduce the notion of twisted symplectic
reflection algebras and describe the category of representations
of such an algebra associated to a non-faithful $G$-action 
in terms of those for faithful actions of $G$.
\end{abstract}

\author{Tatyana Chmutova}
\address{Department of Mathematics, Harvard University,
Cambridge, MA 02138, USA}
\email{chmutova@math.harvard.edu}

\maketitle

\section{Introduction}
The symplectic reflection algebras $\mathcal H_c(G)$ corresponding
to a finite subgroup $G \subset Sp(U)$ were introduced by Etingof
and Ginzburg in \cite{EG}. The representations of these algebras
and of their special cases called rational Cherednik algebras were
studied by many authors (see e.g. \cite{BEG}, \cite{BEG1},
\cite{DO}, \cite{Du}).

In this paper we
\begin{itemize}
    \item   introduce a twisted version of symplectic reflection
    algebras;
    \item consider symplectic reflection algebras corresponding to
    a representation $G \to Sp(U)$ which is not necessarily
    injective.
\end{itemize}

We show that the noninjective case can be reduced to the injective case,
but in this reduction, even if the cocycle we started with was
trivial it might become nontrivial under reduction.

The structure of the paper is as follows. In Section \ref{proj} we recall some 
basic facts about twisted algebras and projective representations of finite groups.
In Section \ref{twrefl} we define twisted reflection algebras. In Section \ref{mod}
we recall some facts about module categories. In Section \ref{reduction} we describe 
the reduction of the noninjective case to injective ones. 

\section{Projective representations of finite groups.}\label{proj}
In this section we  recall some standard facts about projective
representations of finite groups. The main reference for this
section is \cite{K}.

Let $G$ be a finite group and $\psi: G \times G \to \CC^*$ be a
two-cocycle.

\begin{definition}
{\it The twisted group algebra of} $G$ corresponding to 
$\psi$ is an associative algebra
$\CC_{\psi}G$ with a basis $g \in G$ and multiplication given by
\begin{equation}\label{twistproduct}
g_1\circ_{\psi} g_2=\psi(g_1,g_2) g_1 g_2,
\end{equation}
where $g_1g_2$ is the usual product in $G$.
\end{definition}

\begin{definition} An element $g\in G$ is called {\it $\psi$-regular} if
$\psi(g,h)=\psi(h,g)$ for all $h$ in the centralizer of $g$ in
$G$.
\end{definition}

\begin{remark} \label{conjinv}It is known (see \cite{K}, Lemma 3.6.1), that if
$g \in G$ is $\psi$-regular, so is any conjugate of $g$.
\end{remark}

\begin{remark}
Let $\psi$ and $\psi'$ be two cohomologous cocycles. Then any
element $g \in G$ is $\psi$-regular if and only if it is
$\psi'$-regular.
\end{remark}

\medskip

 Let $\tau$ be a representation of $\CC_{\psi}G$. One
can define the character of $\tau$ in the usual way:
$$
\chi_{\tau}(g) =\tr g|_{\tau}.
$$
The characters of projective representations are not necessarily
class functions, i.e. they might be not invariant under
conjugation.
\begin{definition}
A cocycle $\psi$ is called {\it a class-function cocycle} if
$\chi_{\tau}(g)$ is a class function for every representation
$\tau$ of $\CC_{\psi}G$.
\end{definition}

\begin{proposition}[see \cite{K}, Proposition
7.2.2] Every $\CC^*$-valued two cocycle is cohomologous to a
class-function cocycle.
\end{proposition}

Analogously to the case of the trivial cocycle, the values of a
twisted character on the $\psi$-regular elements determine the
isomorphism class of $\tau$.

\begin{proposition}[see \cite{K}, Theorem 7.2.5 (ii)]\label{projchar}
An element $g$ is $\psi$-regular if and only if $\chi_{\tau}(g)
\ne 0$ for some irreducible representation $\tau$ of
$\CC_{\psi}G$.
\end{proposition}

\begin{theorem}[see \cite{K}, Theorem 3.6.3]
Let $\psi$ be a class-function cocycle, and let $\widetilde{C}_1,
\dots, \widetilde{C}_t$ be $\psi$-regular conjugacy classes of
$G$. Let
$$
k_i=\sum_{g \in \widetilde{C}_i} g.
$$
Then $\{k_1,\dots k_t\}$ is a basis of the center of
$\CC_{\psi}G$.
\end{theorem}

\begin{theorem}[see \cite{K}, Theorem 3.6.7]
The number of irreducible representations of $\CC_{\psi}G$ is
equal to the number of $\psi$-regular conjugacy classes of $G$.
\end{theorem}

\begin{example}[see \cite{K}, Theorem 3.7.3]\label{dihedr}
Let $G=I_2(2m)$ be the group of symmetries of $2m$-gon. It is
generated by two elements $s_1$ and $s_2$ with relations
$$
s_i^2=1, \quad (s_1s_2)^{2m}=1.
$$
Denote by $\varepsilon$ a primitive $2m$th root of 1. Every
non-trivial cocycle $\psi': G \times G\to \CC^*$ is cohomologous
to $\psi:G \times G \to \CC^*$ defined by
$$
\psi\left((s_1s_2)^i, (s_1s_2)^j s_1^{\delta}\right)=1, \quad
\psi\left((s_1s_2)^is_1,(s_1s_2)^js_1^{\delta}
\right)=\varepsilon^j
$$
for all $0 \leqslant i,j \leqslant 2m-1$ and $\delta \in\{0,1\}$.

There are $m$ $\psi$-regular conjugacy classes in $G$:
$$
\{s_1s_2, s_2s_1 \}, \quad \{(s_1s_2)^2,(s_2s_1)^2\},\quad \dots,
\quad \{(s_1s_2)^m, (s_2s_1)^m\}.$$
\end{example}

\begin{example}\label{BC} Let $G=S_n \ltimes \left( \ZZ/m\ZZ\right)^n$ be the
generalized symmetric group. It is generated by elements $s_i, ~1
\leqslant i \leqslant n-1$ and $w_j, ~1 \leqslant j \leqslant n$
with relations
$$
s_i^2=w_j^m=1; \quad (s_is_{i+1})^3=1 \text{ for } 1 \leqslant i
\leqslant n-2; \quad (s_is_j)^2=1 \text{ for } |i-j|\geqslant 2;
$$
$$
s_i w_i=w_{i+1}s_i; \quad s_iw_j=w_js_i \text{ for } j \ne i, i+1;
\quad w_iw_j=w_jw_i.
$$

The Schur multiplier of $G$ was computed in \cite{DM} (see also
\cite{MJ}):
$$
H^2(G, \CC^*)=\left\{%
\begin{array}{ll}
    \ZZ/2\ZZ=\{ \gamma \}, & \hbox{if $m$ is odd, $n\geqslant4$;} \\
    \ZZ/2\ZZ \times \ZZ/2\ZZ \times \ZZ/2\ZZ=\{(\gamma, \lambda, \mu )\}, &
    \hbox{if $m$ is even, $n\geqslant 4$;} \\
    \ZZ/2\ZZ \times \ZZ/2\ZZ=\{(\lambda, \mu)\}, &
    \hbox{if $m$ is even, $n=3$;} \\
    \ZZ/2\ZZ=\{(\mu)\}, & \hbox{if $m$ is even, $n=2$;} \\
    \{e\}, & \hbox{otherwise,} \\
\end{array}%
\right.
$$
where $\gamma, \lambda, \mu \in \{\pm 1,1\}.$

Fix an element $(\gamma, \lambda, \mu)\in H^2(G, \CC^*)$. There
exist a cocycle $\psi:G \times G \to \CC^*$ whose class $[\psi]$
in $H^2(G, \CC^*)$ is $(\gamma, \lambda, \mu)$ such that
$\CC_{\psi}G$ is generated by elements $t_i, ~1 \leqslant i
\leqslant n-1$ and $u_j, ~1 \leqslant j \leqslant n$ with
relations
$$
t_i^2=u_j^m=1; \quad (t_it_{i+1})^3=1 \text{ for } 1 \leqslant i
\leqslant n-2; \quad (t_it_j)^2=\gamma \cdot 1 \text{ for }
|i-j|\geqslant 2;
$$
$$
t_i u_i=u_{i+1}t_i; \quad t_iu_j=\lambda u_jt_i \text{ for } j \ne
i, i+1; \quad u_iu_j=\mu u_ju_i \text{ for } i \ne j.
$$

\smallskip

If $[\psi]=(\pm 1, 1,1)$, then the elements $w_n^k$ are
$\psi$-regular for all $1\leqslant k \leqslant m-1$. Indeed, the
centralizer of $g=w_n^k \in G$ is
$$C_G(w^k_j)=S_{n-1} \ltimes (\ZZ/m\ZZ)^{n-1} \times
(\ZZ/m\ZZ)$$ generated by all the $w_i$ and by $s_i, ~1\leqslant i
\leqslant i-2$. From the relations of $\CC_{\psi}G$ it follows
that
$$
\psi(w_i,w_n)=\mu \psi(w_n,w_i)\text{ for } i \ne n \text{ and }$$
$$
\psi(w_n, s_j)=\lambda\psi(s_j,w_n) \text{ for } j \ne n-1.
$$
Hence if $\lambda=\mu=1$, then $\psi(w_n^k,h)=\psi(h, w_n^k)$ for
all $h$ is the centralizer of $w_n^k$.

\smallskip

If $m$ is not divisible by $4$, then the elements $s_i$ are
$\psi$-regular only for $\psi$ cohomologous to the trivial
cocycle. If $m$ is divisible by $4$, then $s_i$ are $\psi$-regular
for $[\psi]=(1,1,\pm 1)$. 

Indeed,
 for odd $m$ the centralizer of $g=s_{n-1}$ is
generated by $s_j$ for $|j-n+1|\geqslant 2$ and by $w_j$ for $j
\ne n-1, n$. From the relations of $\CC_{\psi}G$ it follows that
$$
\psi(s_{n-1}, s_j)=\gamma \psi (s_j, s_{n-1}).
$$
Hence $g$ is regular only for the trivial cocycle.

For even $m$ the centralizer of $g=s_{n-1}$ is generated by $s_j$
for $|j-n+1|\geqslant 2,~$ $w_j$ for $j \ne n-1, n$, and
$(w_{n-1}w_n)^{m/2}s_{n-1}$. The relations for $\CC_{\psi} G$ give
us
$$
\psi(s_{n-1}, s_j)=\gamma\psi(s_j, s_{n-1}), \quad
\psi(s_{n-1},w_j)=\lambda \psi(w_j, s_{n-1}),
$$
$$
\psi(s_{n-1},(w_{n-1}w_n)^{m/2}s_{n-1})
=\mu^{\frac{m^2}{4}}\psi((w_{n-1}w_n)^{m/2}s_{n-1}, s_{n-1})
$$

 If $m$ is not divisible by $4,
~s_{n-1}$ is regular only for trivial cocycle. Otherwise, it is
regular for a cocycle $\psi$ such that its class in $H^2(G,
\CC^*)$ is $[\psi]=(1,1,\pm 1)$.

\medskip

Let $m=2$ in the example above, i.e. let $G=BC_n$:
$$
\begin{picture}(100,15)(0,0)
\put(0,11){\circle*{3}} %
\put(0,11){\line(1,0){15}} %
\put(15,11){\circle*{3}} %
\put(15,11){\line(1,0){15}} %
\put(30,11){\circle*{3}} %
\put(35,11){\dots} %
\put(50,11){\circle*{3}} %
\put(50,11){\line(1,0){15}} %
\put(65,11){\circle*{3}} %
\put(65,12){\line(1,0){15}} %
\put(65,10){\line(1,0){15}} %
\put(80,11){\circle*{3}} %
\put(-1,2){$_{s_1}$} %
\put(14,2){$_{s_2}$} %
\put(29,2){$_{s_3}$} %
\put(58,2){$_{s_{n-1}}$} %
\put(79,2){$_{w_n}$} %
\end{picture}
$$
The conjugacy class of the reflection $w_n$ is $\psi$-regular for
$[\psi]=(\pm 1,1,1)$, all the other reflections are regular only
for the trivial cocycle.
\end{example}

\begin{proposition}\label{regularefl}
Let $G\ne BC_n$ be a finite irreducible Coxeter group. Then all the reflections
in $G$ are regular only for the trivial cocycle.
\end{proposition}

\begin{proof} Consider any two vertices $i$ and $j$ of the Coxeter graph
connected by a simple edge. The corresponding reflections
$s_i$ and $s_j$ generate a subgroup of $G$ isomorphic to $S_3$.
Since $H^2(S_3, \CC^*)=\{e\}$, we can trivialize $\psi$ on this
subgroup.

 Consider any multiple edge
$$
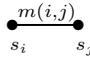
\begin{picture}(30,20)(0,0)
\put(0,8){\circle*{4}} %
\put(0,8){\line(1,0){25}} %
\put(25,8){\circle*{4}} %
\put(-1,0){$_{s_i}$} %
\put(24,0){$_{s_j}$} %
\put(2,13){$_{m(i,j)}$} %
\end{picture}
$$
with $m(i,j)>3$ and $(s_is_j)^{m(i,j)}=1$. In a 
finite irreducible Coxeter group
there are at most $2$ conjugacy classes of reflections: one
corresponding to  short roots and one corresponding to long roots.
So, if there is a $\psi$-regular reflection in $G$, then by Remark
\ref{conjinv} either $s_i$ or $s_j$ is also $\psi$-regular. The
subgroup $G(i,j)$ of $G$ generated by $s_i$ and $s_j$ is
isomorphic $I_2(m(i,j))$. From the Example \ref{dihedr} it follows
that the restriction of $\psi$ on $G(i,j)$ can be trivialized.

One can trivialize $\psi$ on all the edges of the Coxeter graph
simultaneously, i.e. there exist a cocycle $\psi'$ cohomologous to
$\psi$ such that in the twisted group algebra $\CC_{\psi'}G$,
which is generated by elements $t_i$ corresponding to the vertices
of the Coxeter graph of $G$, the following relations are
satisfied:
$$
t_i^2=1, \text{ and } (t_it_j)^{m(i,j)}=1 \text{ for all $i,j$
such that $m(i,j)>2$}.
$$

Let $C$ denote the conjugacy class of $\psi'$-regular reflections
in $G$. Then $\psi'(s_i, s_j)=\psi'(s_j,s_i)$ for all $s_i \in C$
and $s_j$ commuting with $s_i$. This means that $(t_it_j)^2=1$ for
all $(i,j) \in I$, where  $I=\{(i,j)| s_i\in C \text{ and }
m(i,j)=2 \}$.

If $G \ne BC_n$, then $I$ is the set of all pairs of vertices
which are not connected. Since $\CC_{\psi'}G$ and $\CC G$ have
the same generators and relations, they are isomorphic. Hence
$\psi'$ and $\psi$ are cohomologous to the trivial cocycle.
\end{proof}

\section{Twisted symplectic reflection algebras.}\label{twrefl}

Let $G$ be a finite group acting on a symplectic vector space $U$
with symplectic $2$-form $\omega$ and let $\psi:G \times G \to
\CC^*$ be a class-function $2$-cocycle.

\begin{definition} A {\it symplectic reflection} of $U$ is
an element $f \in Sp (U)$ such that $\rank(Id-f)=2$.
\end{definition}

Let $S=S(G,U)\subset G$ be the set of elements of $G$ that act on
$U$ via symplectic reflections and let $S_{\psi}(G,U) \subset
S(G,U)$ be its subset consisting of all $\psi$-regular elements.
For every $s \in S_{\psi}(G,U)$ there is an $\omega$-orthogonal
decomposition $U=\Ker(\Id-s)\oplus \text{Im} (\Id-s)$. Denote by
$\omega_s$ a skew-symmetric form on $U$ which has $\Ker(\Id -s)$
as the radical and coincides with $\omega$ on $\text{Im}(\Id-s)$.

For any conjugation invariant function $c: S_{\psi}(G,U) \to \CC$
one can define the twisted symplectic reflection algebra in the
following way.

\begin{definition} \label{std1}(cf. \cite{EG}, Theorem 1.3)
{\it The twisted symplectic reflection algebra} $\mathcal
H_c(G,U,\psi)$ is the quotient of the semidirect product $(TU)
\rtimes \CC_{\psi}G$ by the two-sided ideal $\mathcal I$ generated
by elements

\begin{equation}
y\otimes x-x\otimes y-\omega(x,y) \cdot 1-  \sum_{s \in
S_{\psi}(G,U)}c(s)\cdot \omega_s(x,y)\cdot s,  \quad \forall x, y
\in U.
\end{equation}
\end{definition}

\begin{remark} In the case of the trivial cocycle $\psi$ and a faithful representation $U$
symplectic reflection algebras were introduced by Etingof and
Ginzburg in \cite{EG}.
\end{remark}

\begin{remark}\label{zero}
If $c(s)=0$ for all $\psi$-regular reflections $s\in S_{\psi}(G,
U)$, then
$$
\mathcal H_c(G,U, \psi)=(SU) \rtimes \CC_{\psi} G,
$$
where $SU$ is the symmetric algebra of $U$.
\end{remark}

\begin{remark}Analogously to the untwisted case one can define
twisted rational Cherednik algebras. 
\end{remark}

\begin{example}

By Proposition \ref{regularefl} and Remark \ref{zero}, the twisted
rational Cherednik algebra $\mathcal H_c(G,U, \psi)$ for $G \ne
BC_n$ is isomorphic to $(SU) \rtimes \CC_{\psi} G$.

 Let $G=BC_n$,
$\psi$ be a nontrivial two-cocycle, and $G'$ be the subgroup of
$G$ generated by the unique $\psi$-regular conjugacy class. Since
the cocycle $\psi$ trivializes on $G'$
$$\mathcal H_c(G,U, \psi)=\CC_{\psi}G \otimes_{\CC
G'}\mathcal H_c(G').$$
\end{example}

\bigskip

The algebra $\mathcal H_c(G, U,\psi)$ has a natural filtration
obtained by placing $\CC_{\psi}(G)$ in grade degree $0$ and $U$ in
grade degree $1$. In the associated graded algebra $\text{\bf gr}
(\mathcal H_c(G,U,\psi))$ any two elements $x,y \in U$ commute.
Hence the tautological imbedding $U \hookrightarrow \text{\bf
gr}(\mathcal H_c(G, U, \psi))$ extends to a surjective graded
algebra homomorphism:
$$
\phi: (SU)\rtimes \CC_{\psi}G \to \text{\bf gr}(\mathcal
H_c(G,U,\psi)).
$$

\begin{theorem}(Poincar\'{e}-Birkhoff-Witt theorem) The homomorphism
$\phi$ defined above is an isomorphism of graded algebras.

\end{theorem}

\begin{proof} Proof is analogous to the untwisted case (see \cite{EG}, Theorem 1.3).
\end{proof}

\begin{definition}Two algebras are called {\it Morita equivalent}
if the categories of modules over these algebras are equivalent.
\end{definition}
In this paper we study the representations of $\mathcal
H_c(G,U,\psi)$ when the action of $G$ on $U$ is not faithful.
Denote by $K$ the kernel of the action. Then we have the following
exact sequence
$$
1\rightarrow K \hookrightarrow G \stackrel{\pi}{\rightarrow} W
\rightarrow 1,
$$
where the action of $W$ on $U$ is faithful. We will show that the
twisted symplectic reflection algebra $\mathcal H_c(G,U,\psi)$
decomposes into a direct sum of subalgebras $\mathcal H_i$ and
each of $\mathcal H_i$ is Morita equivalent to an algebra
$\mathcal H_{c'}(H, U,{\zeta})$ for some subgroup $H \subset W$
and some $2$-cocycle ${\zeta}:H\times H \to \CC^*$.
 To
formulate the statement more precisely we need to introduce the
notion of module categories.

\section{Module categories over $\Rep{W}$}\label{mod}

In what follows we will reproduce the main definitions and facts
about module categories. For more details and proofs see \cite{O}.

\subsection{Module categories and module functors.}

Let $\mathcal C$ be a monoidal category with tensor product
$\otimes: \mathcal C \times \mathcal C \rightarrow \mathcal C$, 
associativity isomorphism $a_{X,Y,Z}$, and functorial isomorphisms
$r_X: X\otimes 1 \to X$ and $l_X: 1\otimes X \to X$ (see \cite{O}
or \cite{BK} for more details).

\begin{definition}(see \cite{BK} Definitions 2.1.1 and 2.1.2)
A monoidal category $\mathcal C$ is called {\it rigid} if every
object in $\mathcal C$ has right and left duals.
\end{definition}

\begin{definition}(cf. \cite{ENO})
{\it A fusion} category $\mathcal C$ is a $\CC$-linear semisimple
rigid monoidal category with finitely many simple objects and
finite dimensional spaces of morphisms, such that the endomorphism
algebra of the neutral object is $\CC$.
\end{definition}

In what follows we assume that $\mathcal C$ is a fusion
category over $\mathcal C$.

\begin{definition}(cf. \cite{O}, Definition 6)
{\it A module category} over $\mathcal C$ is a semisimple category
$\mathcal M$ together with a bifunctor $\otimes : \mathcal C
\times \mathcal M \rightarrow \mathcal M$ and functorial
associativity
 $m_{X,Y,M}:(X \otimes Y)\otimes M
\rightarrow X\otimes (Y \otimes M)$, and unit isomorphisms $l_M:
1\otimes M \rightarrow M$ for any $X,Y \in \mathcal C$ and $M \in
\mathcal M$ such that the following two diagrams commute:

$$
\xymatrix{ (X \otimes 1) \otimes M \ar[dr]_{r_X \otimes id}
\ar[rr]^{m_{X,1,M}}& & X\otimes(1\otimes M) \ar[dl]^{id \otimes l_M} \\
& X \otimes M &}
$$
 and
$$
\xymatrix{ &((X \otimes Y) \otimes Z) \otimes M
\ar[dl]_{a_{X,Y,Z}\otimes id} \ar[dr]^{m_{X\otimes Y,Z,M}}&
\\
(X\otimes (Y \otimes Z))\otimes M \ar[d]^{m_{X,Y\otimes Z,M}} &&
(X \otimes Y)\otimes
(Z\otimes M) \ar[d]_{m_{X,Y,Z\otimes M}} \\
X\otimes((Y\otimes Z) \otimes M) \ar[rr]^{id \otimes m_{Y,Z,M}} &&
X\otimes (Y \otimes( Z\otimes M))}
$$

\end{definition}

\begin{proposition}[\cite{O}, Lemma 2]\label{rigid}
Let $\mathcal M$ be a module category over a fusion category
$\mathcal C$. Then for any $X \in \mathcal C$ and $M_1, M_2 \in
\mathcal M$ we have canonical isomorphisms
$$
\Hom(X \otimes M_1, M_2) \cong \Hom(M_1, ^*X \otimes M_2) \text{
and }
$$
$$
\Hom(M_1, X\otimes M_2) \cong \Hom(X^* \otimes M_1, M_2).
$$
\end{proposition}

\begin{definition}
Let $\mathcal M_1$ and $\mathcal M_2$ be two module categories
over $\mathcal C$. {\it A module functor} from $\mathcal M_1$ to
$\mathcal M_2$ is a functor $F:\mathcal M_1 \rightarrow \mathcal
M_2$ together with functorial isomorphisms $c_{X,M}:F(X\otimes M)
\rightarrow X\otimes F(M)$ for any $X \in \mathcal C$ and $M \in
\mathcal M_1$ such that the following two diagrams commute:
$$
\xymatrix{& F((X\otimes Y)\otimes M) \ar[dr]^{c_{X\otimes Y, M}} \ar[dl]_{Fm_{X,Y,M}}& \\
F(X\otimes(Y \otimes M)) \ar[d]^{c_{X, Y\otimes M}} && (X\otimes
Y)\otimes F(M) \ar[d]_{m_{X,Y,F(M)}}\\
X\otimes F(Y\otimes M) \ar[rr]^{id \otimes c_{Y,M}} && X\otimes
(Y\otimes F(M))}
$$
and
$$
\xymatrix{F(1\otimes M) \ar[rr]^{Fl_M} \ar[dr]_{c_{1,M}}&& F(M)\\
& 1\otimes F(M) \ar[ur]_{l_{F(M)}}}
$$
\end{definition}

\begin{definition}
Two module categories $\mathcal M_1$ and $\mathcal M_2$ over
$\mathcal C$ are {\it equivalent} if there exists a module
functor from $\mathcal M_1$ to $\mathcal M_2$ which is an
equivalence of categories.
\end{definition}

\begin{definition}
Let $\mathcal M_1$ and $\mathcal M_2$ be two module categories
over a fusion category $\mathcal C$. Their {\it direct sum} is the
category $\mathcal M_1 \times \mathcal M_2$ with additive and
module structures defined coordinatewise.
\end{definition}

\begin{definition}
A module category $\mathcal M$ is called {\it indecomposable} if
it is not equivalent to a direct sum of two nontrivial module
categories.
\end{definition}

\subsection{Module categories over $\operatorname{Rep}(W)$.}
In this section we will give a description of all the
indecomposable module categories over the category
$\operatorname{Rep} (W)$ of finite dimensional representations of
a finite group $W$. These module categories are in one-to-one
correspondence with conjugacy classes of pairs $(H, [\zeta])$,
where $H \subset W$ is a subgroup and $[\zeta] \in H^2(H,\CC^*)$ .

Given a pair $(H, \zeta)$ where $\zeta$ is a two-cocycle on $H$,
denote by $\operatorname{Rep}(H,\zeta)$ the category of finite
dimensional representations of $\CC_{\zeta}H$. The usual tensor
multiplication by elements of $\operatorname{Rep}(W)$ defines the
structure of a module category on $\operatorname{Rep}(H,\zeta)$.
Note that if $\zeta$ and $\zeta'$ are cohomologous, then the
corresponding module categories $\operatorname{Rep}(H,\zeta)$ and
$\operatorname{Rep}(H,\zeta')$ are equivalent.

\begin{theorem}[see \cite{O}, Theorem 2] \label{cocycle}
Let $\mathcal M$ be an indecomposable module category over
$\operatorname{Rep}(W)$. Then there exist $H \subset W$ and a
two-cocycle $\zeta:H\times H \to \CC^*$ such that $\mathcal M$ is
equivalent to $\operatorname{Rep}(H,\zeta)$. Two module categories
$\operatorname{Rep}(H_1,\zeta_1)$ and
$\operatorname{Rep}(H_2,\zeta_2)$ are equivalent if and only if
pairs $(H_1, [\zeta_1])$ and $(H_2, [\zeta_2])$ are conjugate
under the adjoint action of $W$.
\end{theorem}

\subsection{Examples.}
1. Consider $W=\ZZ/2\ZZ$. There are two indecomposable module
 categories over $\Rep(\ZZ /2\ZZ)$: $\mathcal M_1$ corresponding
 to the subgroup $e \subset W$ with trivial cocycle and $\mathcal M_2$
 corresponding to the group $W$ itself with trivial cocycle.
 Category $\mathcal M_1$ has one
 irreducible object and category $\mathcal M_2$ has two irreducible objects,
 which are permuted by multiplication by $\sgn \in \Rep(W)$.

2. More generally, let $W=\ZZ /n \ZZ=\langle a\rangle$. Then for
every integer $d$ dividing $n$ there exists a unique indecomposable
module category $\mathcal M_d$ over $\Rep(W)$ with $d$ irreducible
objects. This category corresponds to the subgroup $H \subset W$
of order $d$ with trivial cocycle.

3. Let $G=S_4$ be the symmetric group in 4 elements. The map
$S_4\twoheadrightarrow S_3=W$ defines on $\Rep (S_4)$ a structure of
a module category over $\Rep (S_3)$. The action of $V\in \Rep (S_3)$
is by tensor product. The simple objects of $\Rep (S_4)$ are
$\sf{triv}$, $\sf{sgn}$, the two-dimensional representation $V^2$
and two three-dimensional representations $V_{(1)}^3$ and
$V_{(2)}^3$. The representations $\sf triv$, $\sf sgn$ and $V^2$
generate a module category equivalent to $\Rep(S_3)$, and
$V^3_{(1)}$, $V^3_{(2)}$ generate a module category equivalent to
$\Rep(\ZZ/2\ZZ)$ with the trivial cocycle.

\section{Representations of $\mathcal H_c(G,U,\psi)$.} \label{reduction}

\subsection{Notations}
Let $G$ be a finite group and
$\pi:G{\twoheadrightarrow} W \subset Sp(U)$ be its representation.
For any two-cocycle $\psi: G\times G \to \CC^*$ the category
$\Rep(G,\psi)$ of representations of $\CC_{\psi}G$ is a module
category over $\Rep(W)$ with the action given by tensor product.
Let $\Rep (G,\psi)=\bigoplus_i \mathcal C_i$ be its decomposition
into a direct sum of indecomposable module categories. By Theorem
\ref{cocycle} each of $\mathcal C_i$ is equivalent to a category
$\Rep(H, \zeta)$ for some subgroup $H \subset W$ and some
two-cocycle $\zeta: H \times H \to \CC^*$.

Let $C_1, \dots, C_n$ be conjugacy classes in $W$ such that the
set $C_i \cap H$ is not empty and contains some $\zeta$-regular
elements. Each set $(C_i\cap H)^{(\zeta)}$ consisting of
$\zeta$-regular elements in $C_i\cap H$ is a union of several
conjugacy classes in $H$:
$$(C_i\cap
H)^{(\zeta)}=\bigcup_{j=1}^{m_i}C_i^j.$$

\subsection{Main Lemma.}

\begin{lemma}\label{charcter}
 Let  $\mathcal M$ be a subcategory
 of $\operatorname{Rep}(G, \psi)$
 equivalent to $\Rep(H, \zeta)$ as a module category over $\Rep(W)$. Denote by $F: \mathcal M
 \rightarrow \operatorname{Rep}(H, \zeta)$ the equivalence of these two module categories.
  Fix $g\in G$ and let $\pi(g)$ be its projection on $W$.

\begin{enumerate} \label{mainlemma}
\item[(i)] If $\pi(g) \notin C_i$  for any $i$, then
for every $M\in \mathcal M$
$$\chi_{M}(g)=0.$$

\medskip

\item[(ii)] If $\pi(g) \in C_i$, then for every simple object $M \in \mathcal M$
$$\chi_{M}(g)=\sum_{j=1}^{m_i}\alpha(C_i^j)\chi_{F(M)}(C_i^j)$$
where the coefficients $\alpha(C_i^j)$ do not depend on $M$.
\end{enumerate}
\end{lemma}

\begin{proof}[Proof of the lemma \ref{charcter}]
\noindent
\begin{enumerate}
\item[(i)]
Let $K(\mathcal M)$ be the Grothendieck group of the module
category $\mathcal M $. Denote by $\mathcal{CL}$ the algebra of
class-functions on the set of $\zeta$-regular elements in $H$.
Then we have a map $\Rep(W) \to \mathcal{CL}$, and  by Proposition
\ref{projchar} the action of $\Rep(W)$ in $K(\mathcal M)$ factors
through $\mathcal{CL}$.

Let $f$ be the class function on $W$ which is equal to $1$ at the
conjugacy class of $\pi(g)$ and is zero elsewhere. Then $f$ maps
to $0$ in $\mathcal{CL}$ since $\pi(g)$ is not conjugate to a
$\zeta$-regular element of $H$. Let $A_f \in \Rep(W)$ be the
virtual representation corresponding to $f$. Then $A_f\otimes
V=0$. Taking characters at $g$ we get $f(\pi(g))\chi_M(g)=0$.
Since $f(\pi(g))=1$, we get $\chi_M(g)=0$.

 \medskip

\item[(ii)]
Let $f_i$ be a class function on $W$ which is equal to $1$ on
$C_i$ and is zero elsewhere. Denote by $A_i \in \Rep(W)$ the
corresponding virtual representations. Let $B_1, \dots B_k$ be all
the irreducible representations of $\CC_{\zeta}H$. The action of
$\Rep(W)$ on $\Rep(H, \zeta)$ is given by
$$
A_i \otimes B_j= \sum_{r=1}^k\xi_{i,j}^r B_r.
$$

Consider vectors $v^{(t)} \in \CC^k$ defined by
$v^{(t)}_j=\chi_{B_j}(h_t)$, where $h_t$ are representatives of
$\zeta$-regular conjugacy classes in $H$. These vectors form a
basis in $\CC^k$, since the characters of irreducible
representations form a basis in the space of class functions.
Moreover they are eigenvectors of the matrices
$\Xi_i=\left(\xi_{i,j}^r\right)$:
$$\Xi_i (v^{(t)})=\left\{%
\begin{array}{ll}
    v^{(t)}, & \text{ if }h_t \in C_i; \\
    0, & \text{ otherwise}. \\
\end{array}%
\right.    $$

Let $M_l$ for $1\leqslant l \leqslant k$ be the irreducible
representation of $\CC_{\psi} G$ such that $F(M_l)=B_l$.  Let
$g\in G$ and $\pi(g) \in C_i$. Then the vector $u \in \CC^k$ with
coordinates $u_j=\chi_{M_j}(g)$ satisfies the equations $\Xi_j(u)=
\delta_j^i u$ and hence it is a linear combination of vectors
$v^{(t)}$ such that $h_t$ is in $C_i$.
\end{enumerate}
\end{proof}

\subsection{Main Theorem.}
The category $\Rep(G, \psi)$ decomposes into a direct sum of
indecomposable module categories in the following way:
$$\Rep(G,\psi)=\bigoplus_i \mathcal C_i.$$
In particular, the regular representation $\CC_{\psi}G$ of $G$
decomposes into a direct sum $\CC_{\psi}G=\bigoplus_i
\CC_{\psi}G^{(i)}$, where $\CC_{\psi}G^{(i)}\in \mathcal
Ob(\mathcal C_i)$. Hence $1\in \CC_{\psi}G$ can be written as
$1=\bigoplus_i e_i$, where $e_i \in \CC_{\psi}G^{(i)}$ are
idempotents.

Fix  a conjugation invariant function $c: S_{\psi}(G,U) \to
\CC^*$. Let $\mathcal H_c(G, U, \psi)$ be the corresponding
symplectic reflection algebra. The categories $C_i$ are module
categories over $\Rep(W)$, in particular they are closed under
tensor multiplication by $U$. This means that the constructed
above idempotents $e_i$  are central in $\mathcal H_c(G, U,
\psi)$. Moreover they are minimal with such property since
$\mathcal C_i$ are indecomposable. Hence the algebra $\mathcal
H_c(G, U, \psi)$ can be decomposed in a direct sum of subalgebras
in the following way:
$$\mathcal
H_c(G, U, \psi)=\bigoplus_i \mathcal H_i, $$
where $\mathcal H_i=e_i  \mathcal H_c(G, U, \psi).$

\begin{theorem}\label{equiv}
If the module category $\mathcal C_i$ is equivalent to a category
$\Rep(H, \zeta)$, then the algebra $e_i \mathcal H_c(G, U, \psi)$
is Morita equivalent to the algebra $\mathcal H_{c'}(H, U,
\zeta)$ for   a conjugation invariant function
$c':S_{\zeta}(H,U)\to \CC^*$   defined by

$$
c'(C_i^j)=\frac{\alpha(C_i^j)}{\alpha(e)|C_i^j|}\sum_{\pi
(\widetilde{C})=C_i} c(\widetilde{C})|\widetilde{C}|,
$$
where $C_i^j$ is a conjugacy class of $\zeta$-regular reflections
in $H$ and $|C_i^j|$ and $|\widetilde{C}|$ stands for the number of
elements in the conjugacy class. (See Lemma \ref{mainlemma} for
definition of $\alpha(C_i^j)$ and $\alpha(e)$.)
\end{theorem}

\begin{proof}
Let $\mathcal A$ and $\mathcal B$ be the categories of 
representations of $e_i  \mathcal H_c(G, U, \psi)$ and of
$\mathcal H_{c'}(H, U, \zeta)$ respectively. Denote by
$\overline{F}$ the equivalence $\mathcal C_i \to \Rep(H, \zeta)$.

Let us define a functor $F:\mathcal A \to \mathcal B$ in the
following way. Every $e_i \mathcal H_c(G, U, \psi)$-module $M$ is
in particular a $\CC_{\psi}G$-module which lies in $\mathcal C_i$.
The module $\overline{F}(M)$ has a natural action of $U$ coming
from the corresponding action on $M$ via the functorial
isomorphisms
$$\overline{F}(U\otimes M)\to U\otimes
\overline{F}(M).$$

\begin{claim} The action of $U$
defines  a  structure of $\mathcal H_{c'}(H, U, \zeta)$-module on
$\overline{F}(M)$.
\end{claim}
\begin{proof}
We have to check that if the element
$$y\otimes x-x\otimes
y-\omega(x,y) \cdot 1-  \sum_{s \in S_{\psi}(G,U)}c(s)\cdot
\omega_s(x,y)\cdot s$$ acts by $0$ on $M$ for all  $x, y \in U$,
then
$$y\otimes x-x\otimes y-\omega(x,y) \cdot 1-  \sum_{s \in
S_{\zeta}(H,U)}c'(s)\cdot \omega_s(x,y)\cdot s$$ acts by $0$ on
$\overline{F}(M)$ for all $x, y \in U$.

\bigskip

Let $\displaystyle \rho(x,y)=\omega(x,y) \cdot 1+  \sum_{s \in
S_{\psi}(G,U)}c(s)\cdot \omega_s(x,y)\cdot s$. %
This element can be viewed as  a homomorphisms $\rho$ from $U \otimes U
\otimes M$ to $M$. By Proposition \ref{rigid}, since the category
$\Rep W$ is a fusion category
$$
\Hom(U \otimes U \otimes M , M)\cong \Hom(U\otimes M, U^*\otimes
M).
$$
Denote by $\overline{\rho}$ the map from $U \otimes M $ to $U^*
\otimes M$ corresponding to $\rho$. Let $\overline{\omega}: U^*
\to U$ be the map given by  the symplectic form $\omega$. Then

$$
(\overline{\omega}\otimes \Id)\circ \overline{\rho}=
\Id \otimes \Id -\sum_{s \in S_{\psi}(G,U)}c(s)(\Id \otimes
s-s\otimes s)=
$$
$$
=\Id \otimes \Id - \Id \otimes \left(\sum_{s \in S_{\psi}(G,U)}
c(s)s\right) +\sum_{s \in S_{\psi}(G,U)} c(s) s\otimes s
$$
as a map from $U\otimes M$ to $U\otimes M$.
\bigskip

To prove the claim, it is enough to check that the element
$\displaystyle \sum_{s \in S_{\psi}(G,U)} c(s)s$ acts in any
irreducible representation $\tau$ of $G$  by the same constant as
the element $\displaystyle \sum_{s \in S_{\zeta}(H,U)} c'(s)s$
acts in $\overline{F}(\tau)$.

\bigskip

Let $\tau$ be an irreducible representation of $G$ and
$\widetilde{C}$ be a conjugacy class of $\psi$-regular reflections
in $G$, which projects onto a conjugacy class $C$ in $W$. Then the
element $ \displaystyle \sum_{s\in \widetilde{C}} s$ acts in
$\tau$ by
$$
\lambda=\frac{|\widetilde{C}|\chi_{\tau}(\widetilde{C})}{\dim
\tau},
$$
 where $|\widetilde{C}|$ is the number of elements in
$\widetilde{C}$.

By lemma \ref{mainlemma} $\lambda$ is nonzero only if $C=C_i$ for
some $i$, i.e. if $C\cap H$ contains a $\zeta$-regular element. In
this case
$$
\lambda = \frac{|\widetilde{C}| \displaystyle\sum_{j=1}^{m_i}
\alpha(C_i^j)\chi_{F(\tau)}(C_i^j)}{\alpha(e) \dim F(\tau)}.
$$

Hence the element $\displaystyle \sum_{s \in S_{\psi}(G,U)} c(s)s$
acts in $\tau$ via multiplication by $\Lambda$, where

$$
\Lambda=\sum_{\widetilde{C}}c(\widetilde{C}) \lambda=\sum_i
\sum_{\pi(\widetilde{C})=C_i}c(\widetilde{C})\frac{|\widetilde{C}|
\displaystyle\sum_{j=1}^{m_i}
\alpha(C_i^j)\chi_{F(\tau)}(C_i^j)}{\alpha(e) \dim F(\tau)}=
$$
$$
=\sum_i \sum_{j=1}^{m_i}\frac{|C_i^j|\chi_{F(\tau)}(C_i^j)}{\dim
F(\tau)}\cdot %
\frac{\alpha(C_i^j)}{\alpha(e)|C_i^j|}\sum_{\pi(\widetilde{C})=C_i}
c(\widetilde{C})|\widetilde{C}|.
$$
In the formulas above the sum is taken over those $i$, for which
$C_i$ consist of reflections.

\newpage

Let

$$
c'(C_i^j)=\frac{\alpha(C_i^j)}{\alpha(e)|C_i^j|}\sum_{\pi
(\widetilde{C})=C_i} c(\widetilde{C})|\widetilde{C}|.
$$
Since $\displaystyle \sum_{s\in C_i^j}s$ acts by
$\displaystyle\frac{|C_i^j|\chi_{F(\tau)}(C_i^j)}{\dim F(\tau)}$
in $\overline{F}(\tau)$, the element $\displaystyle\sum_{s \in
S_{\zeta}(H,U)}c'(s)s$ acts in $\overline{F}(\tau)$ by $\Lambda$.

\end{proof}
\begin{claim}
The functor $F$ is an equivalence of categories.
\end{claim}
\begin{proof}
\noindent
\begin{enumerate}

\item  The functor $F$ is surjective on isomorphism
classes of objects in $\mathcal B$. Indeed, any object $N \in
\mathcal B$ is a $\Rep(H, \zeta)$-module, and hence is equal to
$\overline{F}(M)$ for some $\Rep(G, \psi)$-module $M$. Since
$\overline{F}$ is an equivalence of module categories over $\Rep
W$, it respects the $U$-action, so $F(M)=N$ as $\mathcal
B$-module.
\item The map $F: \Hom_{\mathcal A}(M_1,M_2) \to \Hom_{\mathcal B}
(F(M_1),F(M_2))$ is an isomorphism for any $M_1, M_2 \in \mathcal
A$, since $\overline{F}$ is an equivalence of module categories.
\end{enumerate}
\end{proof}
\end{proof}

\end{document}